\documentclass[12pt]{article}
\usepackage{latexsym,amsmath,amssymb}
\usepackage{graphicx}
\newcommand{\fl}{\longrightarrow}
\newfont{\bb}{msbm10 at 12pt}

\def\r{\hbox{\bb R}}
\def\m{\hbox{\bb M}}
\def\h{\hbox{\bb H}}

\def\s{\hbox{\bb S}}
\def\o{\hbox{\bb O}}

\def\gradi{\|\nabla h\|^2}
\def\grad2{\|\nabla h\|^4}

\newcommand{\beq}{\begin{equation}}

\newcommand{\eeq}{\end{equation}}

\newcommand{\mm}{\hbox{\bb{M}}^2(\varepsilon)\times\hbox{\bb{R}}}

\newcommand{\mmm}{\hbox{\bb{M}}^2(\varepsilon)}

\usepackage[latin1]{inputenc}
\usepackage {amsmath}
\usepackage {amsthm}
\usepackage{times}
\usepackage{amscd}
\usepackage{epsf}

\begin{document}

\theoremstyle{plain}\newtheorem{lema}{Lemma}
\theoremstyle{plain}\newtheorem{proposicion}{Proposition}
\theoremstyle{plain}\newtheorem{teorema}{Theorem}
\theoremstyle{plain}\newtheorem{ejemplo}{Example}
\theoremstyle{plain}\newtheorem{observacion}{Remark}
\theoremstyle{plain}\newtheorem{corolario}{Corollary}


\begin{center}
\rule{15cm}{1.5pt} \vspace{0.5cm}

{\Large \bf Complete surfaces of constant curvature\\[2mm] in H$^2\times $R and S$^2\times $R}\\
\vspace{0.5cm} {\large Juan A. Aledo$\,^\dag$, José M. Espinar$\,^\ddag$, José A.
Gálvez$\,^\ddag$}\\
\vspace{0.3cm} \rule{15cm}{1.5pt}
\end{center}

\vspace{.5cm}

\noindent $\mbox{}^\dag$ Departamento de Matem\'{a}ticas, Universidad de Castilla-La
Mancha, 02071 Albacete, Spain; e-mail: juanangel.aledo@uclm.es\vspace{0.2cm}

\noindent $\mbox{}^\ddag$ Departamento de Geometría y Topología, Universidad de
Granada, 18071 Granada, Spain; e-mails: jespinar@ugr.es; jagalvez@ugr.es


\vspace{.3cm}

\begin{abstract}
We study isometric immersions of surfaces of constant curvature into the homogeneous
spaces $\h^2\times\r$ and $\s^2\times\r$. In particular, we prove that there exists
a unique isometric immersion from the standard 2-sphere of constant curvature $c>0$
into $\h^2\times\r$ and a unique one into $\s^2\times\r$ when $c>1$, up to
isometries of the ambient space. Moreover, we show that the hyperbolic plane of
constant curvature $c<-1$ cannot be isometrically immersed into $\h^2\times\r$ or
$\s^2\times\r$.
\end{abstract}

\section{Introduction.}

The Liebmann and Hilbert theorems on surfaces of constant curvature are two of the
most celebrated results in the theory of submanifolds. The first one states that the
round sphere is the unique complete surface of positive constant curvature in the
Euclidean 3-space $\r^3$. The second one asserts that there is no complete surface
of negative constant curvature in $\r^3$.

These results have been extended to the 3-sphere $\s^3$ and the hyperbolic 3-space
$\h^3$ since the Liebmann original proof as well as the Holmgren proof of the
Hilbert theorem also work in these ambient spaces. In these generalizations, it is
basic the fact that the Codazzi equation remains unchanged and the Gauss equation is
very similar to the one of $\r^3$. However, these proofs cannot be extended to other
ambient spaces because the compatibility equations are, in general, completely
different.

Lately, the Hopf theorem on the characterization of the totally umbilical spheres as
the unique topological spheres of constant mean curvature immersed into a
3-dimensional space form, has been extended to the homogeneous spaces $\h^2\times\r$
and $\s^2\times\r$ by Abresch and Rosenberg \cite{AR}. They classify every immersion
of constant mean curvature from a topological sphere into these spaces. The key of
their proof is to show the existence of a holomorphic quadratic form associated with
the Hopf differential for every surface of constant mean curvature.

This result and a lot of other recent works on product spaces have meant a renewed
interest to these 3-manifolds and have turned this research topic into a fashion
one. An up to date reference list of papers on this subject can be found in
\cite{FM}. However, little is known about surfaces of constant (Gaussian) curvature,
even though there are unsolved relevant geometric problems in this context, like the
following one posed by Alencar, do Carmo and Tribuzy \cite{ACT}: which are the
closed immersed surfaces in $\h^2\times\r$ and $\s^2\times\r$ of constant Gaussian
curvature?


We will deal with the previous problem and as a main
goal we will
show that a Liebmann type theorem and a Hilbert type theorem can
be established in $\h^2\times\r$ and $\s^2\times\r$.

The paper is organized as follows. In Section 2 we review the theory of Codazzi
pairs following the ideas of Milnor \cite{Mi} and its relation with holomorphic
quadratic forms.

Section 3 is devoted to the study of the complete revolution surfaces of constant
curvature in $\h^2\times\r$ and $\s^2\times\r$. A local study of these surfaces in
$\h^2\times\r$ can be seen in \cite{MO}. We prove that a complete revolution surface
of constant curvature $K(I)$ in $\h^2\times\r$ must satisfy $K(I)\geq-1$ and find a
parametrization of them. Moreover, we show that, for any positive constant $c$,
there exists a unique complete revolution surface of curvature $K(I)=c$, up to
isometries of $\h^2\times\r$.

For $\s^2\times\r$, we prove the existence and uniqueness, up to isometries of the
target space, of a complete revolution surface of constant curvature $K(I)=c$ for
any $c\geq1$. In addition, the (flat) cylinders are the only complete revolution
surfaces of constant curvature $K(I)<1$ in $\s^2\times\r$.

In Section 4, we prove that for a large class of surfaces of constant curvature
there exists a Codazzi pair related to its induced metric, second fundamental form
and height function. Besides, this pair has constant extrinsic curvature, which
gives us the existence of a holomorphic quadratic form for any surface of positive
constant curvature in $\h^2\times\r$ and constant curvature greater than one in
$\s^2\times\r$.

This fact, as it happens for surfaces of constant mean curvature \cite{AR}, is the
main key to obtain the characterization of the complete surfaces of constant
curvature. Thus, we obtain a Liebmann type theorem, that is, we prove that there
exists a unique complete surface of positive constant curvature in $\h^2\times\r$
and a unique complete surface of positive constant curvature greater than 1 in
$\s^2\times\r$, up to isometries of the ambient space. These complete surfaces are
precisely  the revolution ones given in Section 3.

Finally, in Section 5 we deal with complete surfaces of negative constant curvature,
obtaining the following Hilbert type theorem:  there exists no complete immersion of
constant curvature $K(I)<-1$ into $\h^2\times\r$ or $\s^2\times\r$.

Note that our study shows the existence or non existence of complete surfaces of
constant curvature $K(I)$ in $\h^2\times\r$ for any real number $K(I)$. However, the
cases $-1\leq K(I)< 0$ and $0<K(I)<1$ remain open in $\s^2\times\r$. Anyway, we will
also point out the non existence of complete surfaces of constant curvature in these
cases if some additional assumption is assumed on its height function.


\section{Preliminaries.}

Let us denote by $\r^4_k$, $k=0,1$, the real vector space $\r^4$ endowed with linear
coordinates $(x_1, x_2, x_3,x_4)$ and the metric $\langle,\rangle$ induced by the
quadratic form $\varepsilon\,x_1^2+x_2^2+x_3^2+x_4^2$, where $\varepsilon=1$ (resp.
$\varepsilon=-1$) if $k=0$ (resp. $k=1$).

We will see $\s^2\times\r$ as the submanifold of the Euclidean space $\r^4_0$, given
by
$$
\s^2\times\r=\{(x_1, x_2, x_3,x_4)\in\r^4:\ x_1^2+x_2^2+x_3^2=1\}.
$$

In a similar way, $\h^2\times\r$ will be considered as the submanifold of the
Lorentzian space $\r^4_1$, given by
$$
\h^2\times\r=\{(x_1, x_2, x_3,x_4)\in\r^4:\ -x_1^2+x_2^2+x_3^2=-1,\ x_1>0\}.
$$

Throughout this work we will study complete immersions in $\s^2\times\r$ and
$\h^2\times\r$ of constant Gaussian curvature. In order to do that, we will need
some results about Codazzi pairs and the existence of holomorphic forms. Thus,
following \cite{Mi}, we will remember the necessary ones.

Let $S$ be an orientable surface and $(A,B)$ a pair of real quadratic forms on $S$
such that $A$ is a Riemannian metric. Associated to this pair it is possible to
define its mean curvature and extrinsic curvature, respectively, as
\begin{equation}
\label{CE} H(A,B)=\frac{E g-2F f+G e}{2(EG-F^2)},\qquad
K(A,B)=\frac{eg-f^2}{EG-F^2},
\end{equation}
where
$$
A=E\,du^2+2F\,dudv+G\,dv^2,\qquad B=e\,du^2+2f\,dudv+g\,dv^2,
$$
for local parameters $(u,v)$ on $S$.

In addition, it is said that $(A,B)$ is a Codazzi pair if it satisfies  the
classical Codazzi equations for surfaces in a 3-dimensional space form, that is,
\begin{equation}
\label{ecuacionCodazzi}
e_v-f_u=e\Gamma_{12}^1+f(\Gamma_{12}^2-\Gamma_{11}^1)-g\Gamma_{11}^2,\quad
f_v-g_u=e\Gamma_{22}^1+f(\Gamma_{22}^2-\Gamma_{12}^1)-g\Gamma_{12}^2,
\end{equation}
where $\Gamma_{ij}^k$ are the Christoffel symbols for the Riemannian metric $A$.

It is well known that the Codazzi pair has constant mean curvature if and only if
the $(2,0)$-part of $B$ is a holomorphic 2-form for the conformal structure induced
by $A$ on $S$ \cite[Lemma 6]{Mi}. This fact shows that the Hopf differential is
holomorphic for surfaces of constant mean curvature in 3-dimensional space forms.


It is less known that a Codazzi pair has positive constant extrinsic curvature if
and only if the $(2,0)$-part of $A$ is a holomorphic 2-form for the conformal
structure induced by $B$ on $S$ \cite[Lemma 8]{Mi}. Observe that, in this case,  $B$
or $-B$ is a Riemannian metric. This result was used by Klotz to give a simple proof
of Liebmann's theorem \cite{Mi2}, and it will be our main tool for the
classification of the surfaces with positive constant Gaussian curvature in
$\s^2\times\r$ and $\h^2\times\r$.

\section{Complete revolution surfaces of constant curvature.}\label{seccion2}

In this section we focus our attention on the study of the complete revolution
surfaces of constant curvature in $\mm$, where $\mmm$ denotes $\s^2$ if
$\varepsilon=1$ or $\h^2$ if $\varepsilon=-1$.

It is well known that the special orthogonal group $\s\o(2)$ can be identified with
the subgroup of isometries of $\mm$ (rotations) which preserve the orientation and
fix every point of an axis $\{p\}\times\r$, with $p\in\mmm$.

We observe that if $\varepsilon=-1$ the set of fixed points is exactly
$\{p\}\times\r$. However, if $\varepsilon=1$ the antipodal axis $\{-p\}\times\r$
also remains unchanged under the action of these rotations.

Up to an isometry, we can assume that the axis is given by $\{(1,0,0)\}\times\r$. In
addition, the set ${\cal P}=\{(x_1, x_2, x_3,x_4)\in\mm :\ x_2\geq 0,\,x_3=0\}$
intersects every $\s\o(2)$-orbit once. Hence, we will start with a curve lying in
${\cal P}$ such that it does not cut the axis except at the initial or end points.
Moreover, the curve can only intersect the axis orthogonally; otherwise the
revolution surface would have a pick.

\subsection{Revolution surfaces in H$^2\times$R.}

Let us consider the curve $\alpha(u)=(\cosh k(u),\sinh k(u),0,h(u))\subseteq{\cal
P}$, where $k(u)\geq 0$ and $u$ is the arc length of $\alpha$, that is,
\begin{equation}\label{cero}
k'(u)^2+h'(u)^2=1.
\end{equation}
 Then, its associated revolution surface is given by
$$
\psi(u,v)=(\cosh k(u),\sinh k(u)\cos v,\sinh k(u)\sin v,h(u))
$$
with induced metric
\begin{equation}\label{uno}
 I=du^2+\sinh^2k(u)\,dv^2.
\end{equation}

In general, the Gaussian curvature of a metric $I=du^2+f(u)^2\,dv^2$ is
$K(I)=-f''(u)/f(u)$. Thus, when $K(I)$ is constant, one has
\begin{equation}\label{dos}
\begin{array}
{ll} f(u)=A\,\cos(\sqrt{K(I)}\,u)+B\,\sin(\sqrt{K(I)}\,u)&\mbox{if }K(I)>0,\\[2mm]
f(u)=A\,u+B&\mbox{if }K(I)=0,\\[2mm]
f(u)=A\,\cosh(\sqrt{-K(I)}\,u)+B\,\sinh(\sqrt{-K(I)}\,u)\quad&\mbox{if }K(I)<0,
\end{array}
\end{equation}
for certain real constants $A,B$ which do not vanish simultaneously.

We distinguish three cases in terms of the sign of the constant $K(I)$.
\\[5mm]
\textbf{1.}\  $K(I)>0$:\vspace{2mm}

In this case we have $\sinh k(u)=A\,\cos(\sqrt{K(I)}\,u)+B\,\sin(\sqrt{K(I)}\,u)$,
from (\ref{uno}) and (\ref{dos}).

On the other hand, if the revolution surface is complete then $\psi$ must be, from
the Gauss-Bonnet theorem, the parametrization of a topological sphere, that is,
$\alpha$ is a curve which intersects the axis at its initial and end points.

Let $u_0$ be the initial point. Then, changing $u$ for $u-u_0$, we have $\sinh
k(u)=C\,\sin(\sqrt{K(I)}\,u)$, for a certain real constant $C$.

Thus, one obtains
$$
k(u)={\rm arcsinh}\left(C\,\sin(\sqrt{K(I)}\,u)\right)
$$
and, from (\ref{cero}),
$$
h'(u)^2=1-\frac{C^2\,K(I)\,\cos^2(\sqrt{K(I)}\,u)}{1+C^2\,\sin^2(\sqrt{K(I)}\,u)}.
$$

But $\alpha$ must intersect orthogonally the axis at the initial point, that is,
$h'(0)=0$. Therefore, $C^2\,K(I)=1$ and
$$
h(u)=-\,\sqrt{\frac{1+K(I)}{K(I)}}\ \ {\rm
arctan}\left(\frac{\cos(\sqrt{K(I)}\,u)}{\sqrt{K(I)+\sin^2(\sqrt{K(I)}\,u)}}\right).
$$

The completeness of the immersion follows easily. In particular, we have proven the
existence of a unique complete revolution surface for any positive constant Gaussian
curvature $K(I)$, up to isometries.\\[5mm]
\textbf{2.}\  $K(I)=0$:\vspace{2mm}

From (\ref{uno}) and (\ref{dos}), one has $\sinh k(u)=A\,u+B$. If $\alpha$ does not
cut the axis then the metric of the revolution surface $I=du^2+(A\,u+B)^2\,dv^2$
would only be complete if $\alpha$ is defined for all $u\in\r$. Thus, $A=0$ since
$\alpha$ does not intersects the axis, and $\psi$ is a cylinder around the axis.

On the other hand, if $\alpha$ cuts the axis then, as above, we can assume that it
does at $u=0$ and $\sinh k(u)=A\,u$. Now, using that the intersection is orthogonal
at $u=0$, that is, $h'(0)=0$, and also (\ref{cero})
$$
k(u)={\rm arcsinh}\,u,\qquad h(u)=-1+\sqrt{1+u^2}.
$$

Finally, it is easy to prove that this surface is also complete.
\\[5mm]
\textbf{3.}\  $K(I)<0$:\vspace{2mm}

Now, we have $\sinh k(u)=A\,\cosh(\sqrt{-K(I)}\,u)+B\,\sinh(\sqrt{-K(I)}\,u)$, from
(\ref{uno}) and (\ref{dos}).

First, we consider the case when $\alpha$ cuts the axis. As above, we can assume
that the intersection happens at $u=0$. Then, $\sinh k(u)=C\,\sinh(\sqrt{-K(I)}\,u)$
for a certain real number $C$. Again $h'(0)$ must vanish, and $C^2\,K(I)=-1$, from
(\ref{cero}).

We observe that
$$
k(u)={\rm arcsinh}\left(\frac{1}{\sqrt{-K(I)}}\ \sinh (\sqrt{-K(I)}\,u)\right)
$$
is non negative for $u\geq 0$, and
\begin{equation}\label{kp}
k'(u)=\frac{\sqrt{-K(I)}\ \cosh (\sqrt{-K(I)}\,u)}{\sqrt{-K(I)+\sinh^2
(\sqrt{-K(I)}\,u)}}.
\end{equation}
Since the revolution surface must be complete, the function $k(u)$ is well defined
for all $u\geq0$. Moreover, $k'(u)^2\leq 1$ from (\ref{cero}), and so it follows
that $-K(I)\leq1$ from (\ref{kp}). Hence, $h(u)$ can be calculated as
$$
h(u)=\sqrt{-\frac{1+K(I)}{K(I)}}\ \log\left(\frac{\cosh
(\sqrt{-K(I)}\,u)+\sqrt{-K(I)+\sinh^2 (\sqrt{-K(I)}\,u)}}{1+\sqrt{-K(I)}}\right).
$$

The completeness of these surfaces is a straightforward computation. Thus, we have
proven for every negative constant $K(I)\geq-1$ the existence, up to isometries, of
a unique complete revolution surface of Gaussian curvature $K(I)$ which cuts the
revolution axis.


Finally we consider the case when the curve $\alpha$ does not touch the rotation
axis. Now, the induced metric is given by
$$
I=du^2+\left(A\,\cosh(\sqrt{-K(I)}\,u)+B\,\sinh(\sqrt{-K(I)}\,u)\right)^2\,dv^2
$$
and must be defined for all $u\in\r$ in order to be complete.

If $A^2<B^2$ then the limit of $\sinh
k(u)=A\,\cosh(\sqrt{-K(I)}\,u)+B\,\sinh(\sqrt{-K(I)}\,u)$ is $-\infty$ when $u$
tends to $-\infty$ or $+\infty$, which contradicts that $k(u)\geq 0$. Thus, $A^2$
must be greater than or equal to $B^2$.

If $A^2=B^2$ then, changing $u$ for $-u$ if necessary, we can assume $A=B$. Thus,
$\sinh k(u)=A\,e^{\sqrt{-K(I)}\,u}$ and writing $u$ instead of $u-u_0$ for a
suitable $u_0\in\r$ we have $$\sinh k(u)=e^{\sqrt{-K(I)}\,u}.$$

In addition,
$$
k'(u)=\frac{\sqrt{-K(I)}\,e^{\sqrt{-K(I)}\ u}}{\sqrt{1+e^{2\sqrt{-K(I)}\ u}}}
$$
and $k'(u)^2\leq 1$ from (\ref{cero}). Hence, $K(I)\geq-1$ and the induced metric of
our surface, $I=du^2+e^{2\sqrt{-K(I)}\,u}\,dv^2$, is clearly complete.

If $A^2>B^2$ then $-1<B/A<1$. Thereby, there exist real numbers $\theta$ and $C>0$
such that $A=C\sinh \theta$, $B=C\cosh\theta$, and so, $\sinh
k(u)=C\,\cosh(\sqrt{-K(I)}\ u+\theta)$. Then changing $u$ for
$u+\theta/\sqrt{-K(I)}$ one has
$$
\sinh k(u)=C\,\cosh(\sqrt{-K(I)}\ u)
$$
and
$$
k'(u)=\frac{C\,\sqrt{-K(I)}\,\sinh(\sqrt{-K(I)}\
u)}{\sqrt{1+C^2\,\cosh^2(\sqrt{-K(I)}\ u)}}.
$$
Again the condition $k'(u)^2\leq 1$ is equivalent to $K(I)\geq-1$. Thus, as in the
previous case, the induced metric is complete.

Therefore, we have obtained the existence of complete revolution surfaces of
negative constant Gaussian curvature $K(I)\geq-1$ which do not cut the revolution
axis, in fact, infinitely many for each $K(I)$.



\subsection{Revolution surfaces in S$^2\times$R.}

Let us consider the curve $\alpha(u)=(\sin k(u),\cos k(u),0,h(u))\subseteq{\cal P}$
where $\cos k(u)\geq0$ and $u$ is the arc length of $\alpha$, that is, as above
\begin{equation}\label{cero1}
k'(u)^2+h'(u)^2=1.
\end{equation}

Then, its associated revolution surface is given by
$$
\psi(u,v)=(\sin k(u),\cos k(u)\cos v,\cos k(u)\sin v,h(u))
$$
with induced metric
\begin{equation}\label{uno1}
 I=du^2+\cos^2k(u)\,dv^2.
\end{equation}

We proceed as in the previous case and distinguish three cases in terms of the sign
of the constant $K(I)$.
\\[5mm]
\textbf{1.}\  $K(I)>0$:\vspace{2mm}

Our complete immersion must satisfy $\cos
k(u)=A\,\cos(\sqrt{K(I)}\,u)+B\,\sin(\sqrt{K(I)}\,u)$ from (\ref{dos}) and
(\ref{uno1}). Arguing as in $\h^2\times\r$ we can assume that $\alpha$ intersects
the rotation axis at its initial point, say $u=0$. Thus,
\begin{equation}
\label{seis1} \cos k(u)=C\,\sin(\sqrt{K(I)}\,u)
\end{equation}
for a certain positive real constant $C$.

Moreover, since
$$
k'(u)=-\frac{C\,\sqrt{K(I)}\,\cos(\sqrt{K(I)}\,u)}{\sqrt{1-C^2\,
\sin^2(\sqrt{K(I)}\,u)}}
$$
and $h'(0)=0$, then $C^2\,K(I)=1$.

On the other hand, we know that $\alpha$ cuts the axis (or its antipodal axis) twice
when $K(I)>0$, that is, $\cos k(u)=0$ at the initial and end points of $\alpha$.
Then, from (\ref{seis1}), $u\in [0,\pi/\sqrt{K(I)}]$. In particular, since $\cos
k(u)$ is less than or equal to $1$ in this interval, we have $K(I)\geq 1$.

Therefore, the immersion can be computed with
$$
\begin{array}
{l}{\displaystyle k(u)={\rm arccos}\left(\frac{1}{\sqrt{K(I)}}\ \sin(\sqrt{K(I)}\,u)\right)}\\[2mm]
{\displaystyle h(u)=-\sqrt{\frac{K(I)-1}{K(I)}}\
\log\left(\frac{\cos(\sqrt{K(I)}\,u)+\sqrt{K(I)-\sin^2(\sqrt{K(I)}\,u)}}{1+\sqrt{K(I)}}
\right) }
\end{array}
$$
for $u\in [0,\pi/\sqrt{K(I)}]$.

The completeness of the immersion follows easily. Thus, we have proven the existence
of a unique complete revolution surface for any positive constant Gaussian curvature
$K(I)\geq1$, up to isometries.
%
%
\\[5mm]
\textbf{2.}\  $K(I)=0$:\vspace{2mm}

Now $\cos k(u)=A\,u+B$ from (\ref{dos}) and (\ref{uno1}). If $\alpha$ does not cut
the axis then the metric of the revolution surface $I=du^2+(A\,u+B)^2\,dv^2$ would
only be complete if $\alpha$ is defined for all $u\in\r$, but this is imposible if
$A\neq 0$  because $\cos k(u)=A\,u+B\in[-1,1]$. Thus, $A=0$  and $\psi$ is a
cylinder around the axis.

If $\alpha$ cuts the axis then we can assume that it does at its initial point, say
$u=0$. Then, $\cos k(u)=A\,u$ and $k(u)$ is not well defined for all $u\geq 0$, that
is, the induced metric is not complete.
\\[5mm]
\textbf{3.}\  $K(I)<0$:\vspace{2mm}

In this case $\cos k(u)=A\,\cosh(\sqrt{-K(I)}\,u)+B\,\sinh(\sqrt{-K(I)}\,u)$. So, if
$\alpha$ does not touch the axis and the immersion is complete then $u$ must vary in
$\r$, but this is imposible because $\cos
k(u)=A\,\cosh(\sqrt{-K(I)}\,u)+B\,\sinh(\sqrt{-K(I)}\,u)\not\in[-1,1]$ for all
$u\in\r$.

If $\alpha$ intersects the axis at its initial point $u=0$ then $\cos
k(u)=C\,\sinh(\sqrt{-K(I)}\,u)$ for a certain $C\in\r$. But
$C\,\sinh(\sqrt{-K(I)}\,u)\not\in[-1,1]$ for all $u\geq0$ and the same argument as
in the previous case shows that the immersion cannot be complete.

Therefore, there is no complete revolution surface of constant negative Gaussian
curvature.

\section{A Liebmann type theorem.}

Let $S$ be an orientable surface and $\psi:S\fl\m^2(\varepsilon)\times\r$ an
immersion with unit normal vector field $N$, where $\m^2(\varepsilon)$ denotes
$\s^2$ if $\varepsilon=1$ or $\h^2$ if $\varepsilon=-1$. Let $I=\langle
d\psi,d\psi\rangle$ and $II=\langle -dN,d\psi\rangle$ be the first and second
fundamental forms of $\psi$, and $K(I)$ its Gaussian curvature.

If we denote by $h$ the height function, that is, the fourth coordinate immersion of
$\psi$, and we assume that $K(I)\neq \varepsilon$ at every point on $S$, then one
can consider the new quadratic form
\begin{equation}\label{A}
A=I+\frac{1}{\varepsilon\,K(I)-1}\,dh^2.
\end{equation}
The pair $(A,II)$ satisfies the following essential property
\begin{lema}\label{lema1}
If $A$ is a Riemannian metric on $S$ then the  extrinsic curvature of the pair
$(A,II)$ is given by
$$
K(A,II)=K(I)-\varepsilon.
$$
\end{lema}
\begin{proof}
Let $\eta$ be the unit normal of $\m^2(\varepsilon)\times\r$ in $\r^4_k$ given by
\begin{equation}
\label{eta} \eta(x_1,x_2,x_3,x_4)=(x_1,x_2,x_3,0)
\end{equation}
and $II_{\eta}=\langle -d\eta,d\psi\rangle$ the scalar second fundamental form
induced by $\eta$ on $S$.

Then the Gauss equation asserts that
\begin{equation}
\label{Gauss} K(I)=K(I,II)+\varepsilon K(I,II_{\eta}).
\end{equation}
Let us first compute the extrinsic curvature of the pair $(I,II_{\eta})$.

Let $(u,v)$ be local isothermal parameters on $S$ for the induced metric $I$, that
is, $I=E\,(du^2+dv^2)$, for a positive function $E$. We have
\begin{equation}\label{ss}
\langle -d\eta,d\psi\rangle=\langle -d\psi,d\psi\rangle+dh^2
\end{equation}
and, from (\ref{CE}),
\begin{equation}
\label{CurvaturaEta}
K(I,II_{\eta})=\frac{(h_u^2-E)(h_v^2-E)-(h_u\,h_v)^2}{E^2}=1-\frac{h_u^2+h_v^2}{E}=1-\|
\nabla h\|  ^2,
\end{equation}
 where $\|  \nabla h\| $ is  the modulus of the
gradient of $h$ for the induced metric.

If we take $c=1/(\varepsilon K(I)-1)$ and put $II=e\,du^2+2\,f\,dudv+g\,dv^2$ then,
from (\ref{CE}),  (\ref{Gauss}) and (\ref{CurvaturaEta}), one gets
$$
\begin{array}
{rl}
K(A,II)&{\displaystyle=\frac{e\,g-f^2}{(E+c\,h_u^2)(E+c\,h_v^2)-(c\,h_u\,h_v)^2}=\frac{K(I,II)}{1+c\,\|
\nabla h\| ^2}= }\\[7mm]
&={\displaystyle \frac{K(I,II)\,(\varepsilon
\,K(I)-1)}{\varepsilon\,K(I)-K(I,II_{\eta})}=K(I)-\varepsilon }
\end{array}
$$
as we wanted to prove.
\end{proof}

For the previous isothermal parameters $(u,v)$ the real quadratic form $A$ is given
by
\begin{equation}\label{A2}
A=(E+c\,h_u^2)\,du^2+ 2\,c\,h_u\,h_v\,dudv+(E+c\,h_v^2)\,dv^2,
\end{equation}
where $I=E(du^2+dv^2)$ and $c=1/(\varepsilon K(I)-1)$.

Thus,
$A$ is a Riemannian metric on $S$ if and only if
\begin{equation}
\begin{array}
{c}
0<(E+c\,h_u^2)+(E+c\,h_v^2)=E\,(2+c\,\| \nabla h\| ^2),\\[3mm]
0<(E+c\,h_u^2)(E+c\,h_v^2)-(c\,h_u\,h_v)^2=E^2(1+c\,\| \nabla h\| ^2),
\end{array}
\end{equation}
or equivalently, $1+c\,\| \nabla h\| ^2>0$.

As a consequence we obtain,
\begin{lema}\label{lema2}
Let $S$ be an orientable surface and $\psi:S\fl\m^2(\varepsilon)\times\r$ an
immersion with Gaussian curvature $K(I)$. Then the quadratic form $A$ given by
(\ref{A}) is not a Riemannian metric on $S$ if and only if
\begin{enumerate}
\item[(a)] there exists a point on $S$ satisfying $0\leq K(I)\leq 1$ and $\|\nabla
h\|^2\geq 1-K(I)$, when $\varepsilon=1$,
\item[(b)] there exists a point on $S$ satisfying $-1\leq K(I)\leq 0$ and $\|\nabla
h\|^2\geq 1+K(I)$, when $\varepsilon=-1$,
\end{enumerate}
where $h$ is the height function for the immersion $\psi$.
\end{lema}
\begin{proof}
If we take $e_4=(0,0,0,1)$ then $h=\langle\psi,e_4\rangle$. Thus, for every unit
tangent vector $w$, one has
\begin{equation}\label{gradiente}
\langle\nabla h,w\rangle=\langle d\psi(w),e_4\rangle.
\end{equation}
Therefore  $\| \nabla h\| \leq 1$.

Hence, the condition $1+c\,\| \nabla h\| ^2>0$ is not satisfied if and only if
$c\leq-1$ and $\| \nabla h\| ^2\geq -1/c$, or equivalently, ($a$) and ($b$) are
fulfilled.
\end{proof}

%

Now, we focus our attention on immersions of constant Gaussian curvature. Let us
consider an immersion $\psi:S\fl \m^2(\varepsilon)\times\r$  of constant Gaussian
curvature $K(I)$. For every point $p\in S$ there exists a parametrization of a
neighbourhood of $p$ and a domain $\Omega\subseteq\r^2$ such that the induced metric
is given by
$$
I=\frac{1}{d^2}\,(du^2+dv^2),\qquad d=\frac{1}{2}\,(1+K(I)(u^2+v^2)),
$$
with $(u,v)\in\Omega$.

If we take $II=e\,du^2+2\,f\,dudv+g\,dv^2$ and use (\ref{ss})
\begin{equation}
\label{psiuu}
\begin{array}
{l} {\displaystyle
\psi_{uu}=-\frac{K(I)\,u}{d}\,\psi_u+\frac{K(I)\,v}{d}\,\psi_v+e\,N+\varepsilon\,
(h_u^2-\frac{1}{d^2})\,\eta }\\[4mm]
{\displaystyle
\psi_{uv}=-\frac{K(I)\,v}{d}\,\psi_u-\frac{K(I)\,u}{d}\,\psi_v+f\,N+\varepsilon\,
h_u\,h_v\,\eta }\\[4mm]
{\displaystyle
\psi_{vv}=\frac{K(I)\,u}{d}\,\psi_u-\frac{K(I)\,v}{d}\,\psi_v+g\,N+\varepsilon\,
(h_v^2-\frac{1}{d^2})\,\eta \, .}
\end{array}
\end{equation}

Moreover, if we call $\nu$ to the fourth coordinate immersion of $N$ and bear in
mind that $\langle -dN,\eta\rangle=\langle N,d\eta\rangle=\langle
N,d\psi\rangle-\nu\,dh=-\nu\,dh$, then
\begin{equation}
\label{Nu}
\begin{array}
{l} -N_u=e\,d^2\,\psi_u+f\,d^2\,\psi_v-\varepsilon\,\nu\,h_u\,\eta\\[2mm]
 -N_v=f\,d^2\,\psi_u+g\,d^2\,\psi_v-\varepsilon\,\nu\,h_v\,\eta\\[2mm]
-\eta_u=(h_u^2\,d^2-1)\,\psi_u+h_u\,h_v\,d^2\,\psi_v+\nu\,h_u\,N\\[2mm]
-\eta_v=h_u\,h_v\,d^2\,\psi_u+(h_v^2\,d^2-1)\,\psi_v+\nu\,h_v\,N.
\end{array}
\end{equation}

Thus, the Gauss and Codazzi equations turn into
\begin{equation}
\label{ecuaciones de compatibilidad}
\begin{array}
{lcl} \mbox{Gauss:}&\quad&K(I)=(e\,g-f^2)\,d^4+\varepsilon\,(1-\|\nabla
h\|^2)\\[2mm]
\mbox{Codazzi (1):}&\quad&{\displaystyle e_v-f_u=-\frac{K(I)\,v}{d}\,(e+g)-\varepsilon\,\frac{\nu\,h_v}{d^2} }\\[4mm]
\mbox{Codazzi (2):}&\quad&{\displaystyle
f_v-g_u=\frac{K(I)\,u}{d}\,(e+g)+\varepsilon\,\frac{\nu\,h_u}{d^2} \ .}
\end{array}
\end{equation}

Now, we are in position to prove our main technical result.

\begin{teorema}\label{teorema1}
Let $S$ be an orientable surface and $\psi:S\fl\m^2(\varepsilon)\times\r$ an
immersion of constant Gaussian curvature. If the quadratic form $A$ given by
(\ref{A}) is a Riemannian metric then $(A,II)$ is a Codazzi pair of constant
extrinsic curvature.
\end{teorema}
\begin{proof}
With the previous notation, one has
\begin{equation}
\label{Ayh} A=\left(\frac{1}{d^2}+c\,h_u^2\right)\,du^2+2\,c\,h_u\,h_v\,dudv+\left(
\frac{1}{d^2}+c\,h_v^2\right)\,dv^2
\end{equation}
where $c=1/(\varepsilon\,K(I)-1)$.

On the other hand, the fourth coordinate in (\ref{psiuu}) gives us
\begin{equation}
\label{huu}
\begin{array}
{l} {\displaystyle
h_{uu}=-\frac{K(I)\,u}{d}\,h_u+\frac{K(I)\,v}{d}\,h_v+e\,\nu }\\[4mm]
{\displaystyle
h_{uv}=-\frac{K(I)\,v}{d}\,h_u-\frac{K(I)\,u}{d}\,h_v+f\,\nu }\\[4mm]
{\displaystyle h_{vv}=\frac{K(I)\,u}{d}\,h_u-\frac{K(I)\,v}{d}\,h_v+g\,\nu \, .}
\end{array}
\end{equation}

Then, using (\ref{huu}), a straightforward computation shows that the Christoffel
symbols associated to the Riemannian metric $A$ are given by
\begin{equation}
\label{SimbCristofel}
\begin{array}
{lll} {\displaystyle
\Gamma_{11}^1=\frac{e\,c\,d^2\,\nu\,h_u}{1+c\,\gradi}-\frac{K(I)\,u}{d} }&\quad&
{\displaystyle
\Gamma_{11}^2=\frac{e\,c\,d^2\,\nu\,h_v}{1+c\,\gradi}+\frac{K(I)\,v}{d}
}\\[4mm]
{\displaystyle
\Gamma_{12}^1=\frac{f\,c\,d^2\,\nu\,h_u}{1+c\,\gradi}-\frac{K(I)\,v}{d} }&\qquad&
{\displaystyle
\Gamma_{12}^2=\frac{f\,c\,d^2\,\nu\,h_v}{1+c\,\gradi}-\frac{K(I)\,u}{d}
}\\[4mm]
{\displaystyle
\Gamma_{22}^1=\frac{g\,c\,d^2\,\nu\,h_u}{1+c\,\gradi}+\frac{K(I)\,u}{d} }&\quad&
{\displaystyle
\Gamma_{22}^2=\frac{g\,c\,d^2\,\nu\,h_v}{1+c\,\gradi}-\frac{K(I)\,v}{d} \,. }
\end{array}
\end{equation}

Therefore, from the Gauss equation (\ref{ecuaciones de compatibilidad})
$$
\begin{array}
{rcl} e\,\Gamma_{12}^1+f\,(\Gamma_{12}^2-\Gamma_{11}^1)-g\,\Gamma_{11}^2&=&
{\displaystyle
-\frac{K(I)\,v}{d}\,(e+g)-\frac{(e\,g-f^2)\,c\,d^2\,\nu\,h_v}{1+c\,\gradi} }\\[4mm]
&=& {\displaystyle -\frac{K(I)\,v}{d}\,(e+g)-\varepsilon\,\frac{\nu\,h_v}{d^2} .}
\end{array}
$$

The first Codazzi equation (\ref{ecuaciones de compatibilidad}) asserts that the
first equation in (\ref{ecuacionCodazzi}) is satisfied. Analogously, the second
equation in (\ref{ecuacionCodazzi}) is also satisfied and $(A,II) $ is a Codazzi
pair. Finally, from Lemma \ref{lema1}, this pair has constant extrinsic curvature.
\end{proof}

It should be observed that the existence of this Codazzi pair in
$\m^2(\varepsilon)\times\r$, for immersions of constant Gaussian curvature, does not
only depend on the Codazzi equations as it happens in the 3-dimensional space forms
but it also depends heavily on the Gauss equation.

For an immersion of positive constant Gaussian curvature $K(I)$ satisfying
$K(I)-\varepsilon>0$, Lemma \ref{lema2} states that  the quadratic form $A$ is
Riemannian. In addition, from Theorem \ref{teorema1}, $(A,II)$ is a Codazzi pair of
positive constant extrinsic curvature, which  assures  the existence of a
holomorphic quadratic form \cite[Lemma 8]{Mi}.

\begin{corolario}\label{corolario1}
Let $S$ be an orientable surface and $\psi:S\fl\m^2(\varepsilon)\times\r$ an
immersion of positive constant Gaussian curvature $K(I)$ such that
$K(I)-\varepsilon>0$. If we consider $S$ as the Riemann surface with the conformal
structure induced by its second fundamental form then
$$
Q=\left(\langle \psi_z,\psi_z\rangle+\frac{1}{\varepsilon
K(I)-1}\,h_z^2\,\right)\,dz^2
$$
is a holomorphic quadratic form, where $z$ denotes a local conformal parameter on
$S$.
\end{corolario}

The existence of this holomorphic quadratic form is the main key for the
classification of the immersions of positive constant curvature.

\begin{teorema}
Given a real constant $K(I)$, there exists, up to isometries, a unique complete
surface of constant Gaussian curvature $K(I)>1$ in $\s^2\times\r$ and a unique
complete surface of constant Gaussian curvature $K(I)>0$ in $\h^2\times\r$.

In addition, these surfaces are rotationally symmetric.
\end{teorema}
\begin{proof}
The existence part of this result is showed in Section \ref{seccion2}. Therefore, we
only need to prove uniqueness.

Let $S$ be a surface and $\psi:S\fl\m^2(\varepsilon)\times\r$ a complete immersion
with positive constant Gaussian curvature $K(I)$, satisfying $K(I)>1$ for
$\varepsilon=1$. Since $K(I)$ is positive, $S$ is compact. Thus, as every compact
surface immersed in $\m^2(\varepsilon)\times\r$ is orientable, then $S$ is a
topological sphere from the Gauss-Bonnet theorem.

On the other hand, a holomorphic quadratic form on a topological sphere must vanish
identically. Therefore, from Corollary \ref{corolario1}, $Q\equiv 0$, that is, the
(2,0)-part of $A$ for the conformal structure induced by $II$ vanishes and $A$, $II$
are conformal.

Thus, there exists a function $\lambda$ such that $II=\lambda\, A$. Hence, the
extrinsic curvature of the pair $(A,II)$ is given by $K(A,II)=\lambda^2$ and Lemma
\ref{lema1} implies that $\lambda^2=K(I)-\varepsilon$.

We can assume that $S$ is the Riemann sphere $\r^2\cup\{\infty\}$ with induced
metric
$$
I=\frac{1}{d^2}\,(du^2+dv^2),\qquad d=\frac{1}{2}\,(1+K(I)(u^2+v^2)),\qquad
(u,v)\in\r^2.
$$
In addition, changing the sign of the unit normal $N$ if necessary, we can also
assume that $II$ is positive definite, that is, $\lambda=\sqrt{K(I)-\varepsilon}$.

Now, if we examine the fourth coordinate in (\ref{Nu})
\begin{equation}\label{nugrad}
\frac{-\nu_u}{1+c\,(1-\nu^2)}=\lambda\,h_u,\qquad\frac{-\nu_v}{1+c\,(1-\nu^2)}=
\lambda\,h_v,\qquad\nu^2=1-\gradi,
\end{equation}
where $c=1/(\varepsilon K(I)-1)$.

Since $S$ is compact, there exists a point $p$ on $S$ which is a minimum for the
height function. Hence, up to an isometry, we can suppose $p=(0,0)$, $h(0,0)=0$,
$h_u(0,0)=0=h_v(0,0)$ and, from (\ref{nugrad}), $\nu(0,0)=1$ because $II$ is
positive definite.

We observe that $c\in(0,\infty)$ if $\varepsilon=1$ and $c\in(-1,0)$ if
$\varepsilon=-1$. So, by integrating the two first equations in (\ref{nugrad}), it
is easy to obtain
\begin{equation}\label{nunu}
\begin{array}
{lll}
{\displaystyle\nu=\sqrt{K(I)}\,\tanh\left(-\sqrt{\frac{K(I)}{K(I)-1}}\,h+c_1\right)
}&\qquad&\mbox{if }\varepsilon=1,\\[3mm]
{\displaystyle\nu=\sqrt{K(I)}\,\tan\left(-\sqrt{\frac{K(I)}{1+K(I)}}\,h+c_2\right)
}&\qquad&\mbox{if }\varepsilon=-1,
\end{array}
\end{equation}
where $c_1$ and $c_2$ are constant which can be computed using that $h(0,0)=0$ and
$\nu(0,0)=1$.

Therefore, (\ref{huu}) turns into
\begin{equation}
\label{hanalitica}
\begin{array}
{l} {\displaystyle
h_{uu}=-\frac{K(I)\,u}{d}\,h_u+\frac{K(I)\,v}{d}\,h_v+\sqrt{K(I)-\varepsilon}\,\left(\frac{1}{d^2}+c\,h_u^2\right)\,\nu }\\[4mm]
{\displaystyle
h_{uv}=-\frac{K(I)\,v}{d}\,h_u-\frac{K(I)\,u}{d}\,h_v+c\,\sqrt{K(I)-\varepsilon}\,h_u\,h_v\,\nu }\\[4mm]
{\displaystyle
h_{vv}=\frac{K(I)\,u}{d}\,h_u-\frac{K(I)\,v}{d}\,h_v+\sqrt{K(I)-\varepsilon}\,\left(
\frac{1}{d^2}+c\,h_v^2\right)\,\nu \, .}
\end{array}
\end{equation}
where we have used (\ref{Ayh}), $II=\lambda\,A$ and that $\nu$ is given by
(\ref{nunu}).

Now, we remark that (\ref{hanalitica}) ensures that $h$ is analytic from Bernstein's
analyticity theorem \cite{Mo} and every derivative of $h$ at a point $q=(u_0,v_0)$
can be computed using the value of $h$ and its first derivatives at $q$. But, since
$h(0,0)=h_u(0,0)=h_v(0,0)=0$, then the height function of our immersion $h$ is
uniquely determined by analyticity.

To obtain uniqueness of our immersion we need to analyze the three first coordinates
of $\psi$. This is equivalent to know the normal vector field $\eta$. Thus, let us
assume the existence of two complete immersions $\psi_1,\psi_2:S\fl
\m^2(\varepsilon)\times\r$ with the same positive constant Gaussian curvature and
height function such that the normals given by (\ref{eta}) are respectively
$\eta_1,\,\eta_2$.

By ignoring the fourth coordinate, $\eta_1$ and $\eta_2$ can be considered as maps
into the Riemannian surface of constant curvature $\m^2(\varepsilon)$. Moreover, the
induced metrics $\langle d\eta_1,d\eta_1\rangle$ and $\langle
d\eta_2,d\eta_2\rangle$ agree at every point, from (\ref{Nu}) and (\ref{nunu}), and
 coincide with the induced metric at $(0,0)$.

Therefore, $\eta_i$ can be seen as an immersion from a connected neighbourhood
$\cal{U}$ of $(0,0)$, $i=1,2$, and $\eta_2 \circ \eta_1^{-1} :\eta_1({\cal
U})\subseteq \m^2(\varepsilon)\fl\m^2(\varepsilon)$ is an isometry. That is,
$\eta_1$ and $\eta_2$ agree on $\cal{U}$ up to an isometry of $\m^2(\varepsilon)$.

Observe that the three first coordinates in (\ref{psiuu}) show the analyticity of
$\eta$, see \cite{Mo}, since $N$ is analytically computed from $\eta,\ \psi_u$ and
$\psi_v$. Hence, $\eta_1$ and $\eta_2$ agree not only locally on ${\cal U}$ but
globally up to an isometry of $\m^2(\varepsilon)$. Thus, $\psi_1$ and $\psi_2$
coincide up to an isometry of $\m^2(\varepsilon)\times\r$.
%
%
%
\end{proof}
This result classifies all complete immersions of positive constant curvature into
$\h^2\times\r$. The same can be said for immersions of constant curvature greater
than $1$ into $\s^2\times\r$.

This condition on the curvature agrees with the one
for
immersions of positive constant curvature into the 3-dimensional space forms $\h^3$
and $\s^3$. That is, the immersions of positive constant curvature into $\h^3$ and
positive constant curvature greater than $1$ into $\s^3$ can be jointly studied.
Observe that these are the cases of positive curvature in space forms where the
immersions are of elliptic type. However, the study of complete immersions of
positive constant curvature less than or equal to one into $\s^3$ is very different.

In our case, it can be said that there do not exist complete surfaces of positive
constant curvature $K(I)<1$ if a certain constraint about the gradient of the height
function is assumed.

\begin{lema}
Let $S$ be a compact surface. There does not exist any immersion
$\psi:S\fl\s^2\times\r$ of positive Gaussian curvature $K(I)$ satisfying $\|\nabla
h\|^2 < 1-K(I)$, where $h$ is the height function.

In particular, there does not exist any complete immersion  into $\s^2\times\r$ of
positive constant Gaussian curvature satisfying $\|\nabla h\|^2 < 1-K(I)$.
\end{lema}
\begin{proof}
Arguing as in the above theorem, $S$ must be a topological sphere. On the other
hand, from Equations (\ref{Gauss}) and (\ref{CurvaturaEta}), one has $K(I,II)<0$,
that is, $II$ is a Lorentzian metric. But this is a contradiction since there do not
exist Lorentzian metrics on a sphere.
\end{proof}

We also know that every slab $\s^2\times\{t_0\}$ of $\s^2\times\r$ is a surface of
constant curvature one. Yet, it is not known the existence of other complete
immersions of constant curvature $K(I)\in(0,1]$. In our opinion, it would be very
interesting to investigate if the slabs are the unique complete immersions of
positive constant curvature $K(I)\leq1$ into $\s^2\times\r$.

There are also two natural questions regarding the quadratic form $Q$ in Corollary
\ref{corolario1} that seem interesting. One is whether the surfaces of constant
curvature are the only surfaces in $\h^2\times\r$ and $\s^2\times\r$ with
holomorphic $Q$. This was proved to be false for constant mean curvature surfaces
and the Abresch-Rosenberg holomorphic quadratic form in \cite{FM2}. The second one,
inspired by \cite{FM}, is whether $Q$ comes from a geometrically defined harmonic
Gauss map on a surface of constant curvature in $\h^2\times\r$ or $\s^2\times\r$, at
least for some special value of the Gaussian curvature $K(I)$.

\section{A Hilbert type theorem.}

We devote this section to the study of complete surfaces with negative constant
curvature. As in the classical Hilbert theorem, we will show the non existence of
complete immersions with a certain negative constant curvature.

\begin{teorema}\label{teorema3}
There is no complete immersion of constant Gaussian curvature $K(I)<-1$ into
$\h^2\times\r$ or $\s^2\times \r$.
\end{teorema}
\begin{proof}
Let $S$ be a surface and $\psi:S\fl\m^2(\varepsilon)\times\r$ an immersion of
constant Gaussian curvature $K(I)<-1$. Let us assume that $\psi$ is complete.

From Lemma \ref{lema2} the quadratic form $A$ given by (\ref{A}) is a Riemannian
metric. In addition, we can see that $A$ is a complete metric. Let us consider a
unit vector $w$ and call $c=1/(\varepsilon K(I)-1)$
\begin{itemize}
\item if $\varepsilon=1$ then $-1/2<c<0$ and $A(w)=\langle w,w\rangle+c\,\langle \nabla h,w\rangle^2\geq
(1+c)\,\langle w,w\rangle$ since $\langle \nabla h,w\rangle^2\leq\langle
w,w\rangle$,
\item if $\varepsilon=-1$ then $c>0$ and $I\leq A$.
\end{itemize}
That is, $A$ is greater than or equal to a complete metric. Hence, $A$ is complete.

Now, we compute the Gaussian curvature $K(A)$ of the metric $A$. If we write
$A=E\,du^2+2\,F\,dudv+G\,dv^2$ then it is easy to see that
$$
\begin{array}
{lll}
(E\,G-F^2)\,K(A)&=&E\,\left(\,(\Gamma_{22}^1)_u-(\Gamma_{12}^1)_v+\Gamma_{22}^1\,\Gamma_{11}^1+\Gamma_{22}^2\,
\Gamma_{12}^1-(\Gamma_{12}^1)^2-\Gamma_{12}^2\,\Gamma_{22}^1\,\right)+\\[2mm]
& &F\left(\,(\Gamma_{22}^2)_u-(\Gamma_{12}^2)_v+
\Gamma_{11}^2\,\Gamma_{22}^1+\Gamma_{22}^2\,\Gamma_{12}^2-\Gamma_{12}^1\,
\Gamma_{12}^2-\Gamma_{12}^2\, \Gamma_{22}^2\,\right),
\end{array}
$$
where $\Gamma_{ij}^k$ are the Christoffel symbols of $A$.

Then, from (\ref{Ayh}), (\ref{huu}) and (\ref{SimbCristofel}), a straightforward
computation gives us
$$
\begin{array}
{lll} K(A)&=&{\displaystyle
\frac{K(I)\,(1+c\,h_u\,d^2\,(h_u+u\,d\,(e+g)\,\nu))-c\,h_u\,d^4\,(h_u\,d^2(e
\,g-f^2)+(f_v-g_u)\,\nu)}{1+c\,\gradi} }\\[3mm]
 & &{\displaystyle +\,\frac{c\,d^4\,(e\,g-f^2)\,\nu^2}{(1+c\,\gradi)^2} }
\end{array}
$$
where we have used $\gradi=d^2\,(h_u^2+h_v^2)$.

From the Gauss  and second Codazzi equations (\ref{ecuaciones de compatibilidad})
and the third equation in (\ref{nugrad}), one has
$$
K(A)=\frac{(1+c)\,K(I)-\varepsilon\,c\,(1-\gradi)^2}{(1+c\,\gradi)^2}.
$$

Moreover, if we consider $K(A)$ as a function of $\gradi$ then $K(A)$ is monotonous
and evaluating at 0 and 1, one has
\begin{equation}
\label{K(A)} K(I)-1\leq K(A)\leq K(I)+1<0.
\end{equation}

But, since $(A,II)$ is a Codazzi pair with negative constant extrinsic curvature,
from Theorem \ref{teorema1}, and $A$ is complete, then the infimum of $|K(A)|$ is
zero on $S$ (see \cite[p. 172]{We} and \cite{Wi}) which contradicts (\ref{K(A)}).
Thus, the immersion $\psi$ of constant negative Gaussian curvature $K(I)<-1$ cannot
be complete as we wanted to show.
\end{proof}

As we have proven in Section \ref{seccion2} there exist complete immersions of every
constant curvature $K(I)\geq-1$ into $\h^2\times\r$. Thus, this result cannot be
improved for surfaces of constant curvature in $\h^2\times\r$.

However, we do not know the existence of complete surfaces of constant curvature
satisfying $-1\leq K(I)<0$ in $\s^2\times\r$. In this case, we can also observe that
there is no complete surface if a certain constraint about the gradient of the
height function is assumed.

\begin{lema}
There is no complete immersion of constant Gaussian curvature $-1\leq K(I)<0$ into
$\s^2\times\r$ satisfying
$$
\gradi\leq c_0<1+K(I)\qquad \mbox{or }\qquad \gradi\geq c_0>1+K(I),
$$
for a constant $c_0\in\r$.
\end{lema}
\begin{proof}
It can be shown as in Theorem 3 that $A$ is a complete Riemannian metric. So, as
above, the infimum of $|K(A)|$ must be zero. But $K(A)$ only vanishes at $\gradi
=1+K(I)$. Hence, using the bound of $\gradi$ we have that $|K(A)|$ is greater than a
positive constant, which is a contradiction.
\end{proof}

In our opinion, it would be very interesting to investigate the existence of
complete surfaces of constant curvature $-1\leq K(I)<0$ in $\s^2\times\r$. For these
surfaces the metric $A$ is complete and $(A,II)$ is a Codazzi pair of constant
negative extrinsic curvature, which guarantees the existence of a Tschebysheff net
(see \cite[p. 172]{We}) on $\r^2$ (the universal cover of the surface $S$). Observe
that the existence of a Tschebysheff net has been of a great interest for the study
of some families of surfaces of constant curvature as flat surfaces in $\s^3$ or
pseudo-spherical surfaces in the Euclidean 3-space and its relation with the theory
of integrable systems.

Finally, we remark that there exist infinitely many complete flat surfaces in
$\m^2(\varepsilon)\times\r$. Given a regular curve $\alpha(t)$ in
$\m^2(\varepsilon)$ parametrized by the arc length such that $t\in\r$ then
$\psi(t,s)=(\alpha(t),s)$, $(t,s)\in\r^2$, is a complete flat immersion into
$\m^2(\varepsilon)\times\r$. These surfaces are cylinders on the curve $\alpha$.

But every complete flat surface is not necessarily a cylinder. We have shown in
Section \ref{seccion2} that there also exists a revolution complete flat immersion
which is not a cylinder into $\h^2\times\r$. In fact, it is proven in \cite{GM} the
existence of infinitely many complete flat surfaces in $\h^2\times\r$ which can be
considered as a graph on $\h^2$.

\footnotesize J.A. Aledo was partially supported by  Ministerio de Education y
Ciencia Grant No. MTM2004-02746 and Junta de Comunidades de Castilla-La Mancha,
Grant no. PAI-05-034.

J.M. Espinar and J.A. Gálvez were partially supported by Ministerio de Education y
Ciencia Grant No. MTM2004-02746 and Junta de Andalucía Grant No. FQM325.

\end{document}